\documentclass{amsart}
\usepackage{amssymb}

\newtheorem{theorem}{Theorem}

\newtheorem{lemma}{Lemma}
\newtheorem{proposition}{Proposition}

\newcommand{\bee}[1]{\begin{equation}\label{#1}}
\newcommand{\beq}[1]{\begin{eqnarray}\label{#1}}
\newcommand{\ene}{\end{equation}}
\newcommand{\eqe}{\end{eqnarray}}

\begin{document}
\title{On identities of infinite dimensional Lie superalgebras}

\author[D. Repov\v s, and M. Zaicev]
{Du\v san Repov\v s and Mikhail Zaicev}

\address{Du\v san Repov\v s \\Faculty of Mathematics and Physics, and
Faculty of Education, University of Ljubljana,
P.~O.~B. 2964, Ljubljana, 1001, Slovenia}
\email{dusan.repovs@guest.arnes.si}

\address{Mikhail Zaicev \\Department of Algebra\\ Faculty of Mathematics and
Mechanics\\  Moscow State University \\ Moscow,119992 Russia}
\email{zaicevmv@mail.ru}

\thanks{The first author was supported by the Slovenian Research Agency
grants
P1-0292-0101 and J1-4144-0101.
The second author was partially supported by RFBR grant No 13-01-00234a. We thank the referee for several comments and suggestions.}



\begin{abstract}
We study codimension growth of infinite dimensional Lie superalgebras over an
algebraically closed field of characteristic zero. We prove that if a Lie superalgebra 
$L$ is a Grassmann envelope of a finite dimensional simple Lie algebra then the 
PI-exponent of $L$ exists and it is a positive integer.
\end{abstract}


\maketitle

\section{Introduction}\label{intr}

We shall
consider algebras over a field $F$ of characteristic zero. One of the approaches
in the investigations
of associative and non-associative algebras is to study numerical
invariants associated with their identical relations. Given an algebra $A$, we can
associate the sequence of its codimensions $\{c_n(A)\}_{n\in \mathbb N}$ 
(all
notions and definitions will be given in the next section).

This sequence gives some information not only about identities of $A$ but 
also about structure of $A$. For example, $A$ is nilpotent if and only if 
$c_n(A )=0$ for all large enough $n$. If $A$ is an associative non-nilpotent
$F$-algebra then $A$ is
commutative if and only if $c_n(A)=1$ for all $n\ge 1$.

For an associative algebra $A$  with a non-trivial polynomial identity the sequence $c_n(A)$ is exponentially bounded by the celebrated Regev's Theorem \cite{Reg72}
while $c_n(A)=n!$ if $A$ does not satisfy any non-trivial polynomial identity.
In the
non-associative case the sequence of codimensions may
have even 
faster growth. 
For example, if $A$ is an 
absolutely free algebra then
$$
c_n(A)=a_n n!
$$
where
$$
a_n=\frac{1}{2}{2n-2\choose n-1}
$$
is the Catalan number, i.e. the number of all possible arrangements of brackets in 
the word of length $n$.

For Lie algebra $L$ the sequence $\{c_n(L)\}_{n\in \mathbb N}$ is not exponentially bounded in
general even if $L$ satisfies non-trivial Lie identities (see for example \cite{P}). 
Nevertheless, a class of Lie algebras with exponentially bounded codimensions is sufficiently wide.
It includes in particular, all finite dimensional algebras
\cite{B-Dr,GZ-TAMS2010}, Kac-Moody algebras \cite{Z1,Z2}, infinite
dimensional simple Lie algebras of Cartan type \cite{M1}, Virasoro algebra, and many others.

In the case when
$\{c_n(A)\}_{n\in \mathbb N}$ is exponentially bounded,
the upper and the lower limits of
the sequence $\{\sqrt[n]{c_n(A)}\}_{n\in \mathbb N}$ exist and
a natural question arises: does the ordinary limit
$$
\lim_{n\to\infty} \sqrt[n]{c_n(A)}
$$
exist? 
In case of existence we call this limit $exp(A)$ or PI-exponent of $A$.

 Amitsur conjectured in the
1980's
that for any associative P.I. algebra such a limit exists and 
it is a non-negative integer. This conjecture was confirmed first for verbally
prime P.I. algebras in \cite{BR1,R2},  and later in the general case in \cite{GZ1,GZ2}. 
For Lie algebras a series of positive results was obtained for finite dimensional 
algebras \cite{GRZ2,GRZ1,Z3}, 
for algebras with nilpotent commutator subalgebras \cite{PM}, 
for affine Kac-Moody algebras \cite{Z1,Z2}, and some other classes (see \cite{M2}). 
For Lie superalgebras 
there exist
only partial results  \cite{MZ4,MZ3,MZ2,MZ1}.

On the other hand it was shown in \cite{ZM} that there exists a Lie algebra $L$ with
$$
3.1 < \liminf_{n\to\infty} \sqrt[n]{c_n(L)} \le
\limsup_{n\to\infty} \sqrt[n]{c_n(L)} < 3.9 \quad .
$$
This algebra $L$ is soluble and almost nilpotent, i.e. it contains a nilpotent ideal
of finite codimension. In the
general non-associative case there exists,
for any real number $\alpha >1$,
 an algebra $A_\alpha$ such that
$$
\lim_{n\to\infty} \sqrt[n]{c_n(A_\alpha)}=\alpha.
$$
(see \cite{GMZ}). Note also that by
a recent result \cite{GZ} there exist finite
dimensional Lie superalgebras with a
fractional limit $\sqrt[n]{c_n(L)}$.

In the
present paper we
shall
study Grassmann envelopes of finite dimensional simple Lie
algebras. Our main result is the following theorem:
\medskip
\begin{theorem}\label{t1}
Let $L_0\oplus L_1$ be a finite dimensional simple Lie algebra over an
algebraically closed field $F$ of characteristic zero with some ${\mathbb Z}_2$-grading. Let also $\widetilde L=L_0\otimes G_0\oplus L_1\otimes G_1$ be the
Grassmann envelope of $L$. Then the limit
$$
exp(\widetilde L)=\lim_{n\to \infty}\sqrt[n]{c_n(\widetilde L)}
$$
exists and 
is a positive integer. Moreover, $exp(\widetilde L)=\dim L$.
\end{theorem}

\medskip

Another result of our paper concerns  graded identities. Since any Lie superalgebra
$L$ is ${\mathbb Z}_2$-graded one can consider ${\mathbb Z}_2$-graded identities of
$L$ and the corresponding graded codimensions $c_n^{gr}(L)$. We shall 
prove that graded codimensions have  similar properties.
\medskip

\begin{theorem}\label{t2}
Let $L=L_0\oplus L_1$ be a finite dimensional simple Lie algebra over an
algebraically closed field $F$ of characteristic zero with some ${\mathbb Z}_2$-grading. Let also $\widetilde L=L_0\otimes G_0\oplus L_1\otimes G_1$ be a
Grassmann envelope of $L$. Then the limit
$$
exp^{gr}(\widetilde L)=\lim_{n\to \infty}\sqrt[n]{c_n^{gr}(\widetilde L)}
$$
exists and 
is a non-negative integer. Moreover, $exp^{gr}(\widetilde L)=\dim L$.
\end{theorem}
\medskip

In other words, both PI-exponent $exp(\widetilde L)$ and graded PI-exponent
$exp^{gr}(\widetilde L)$ exist, they are integers and they coincide. Note that for an
arbitrary
${\mathbb Z}_2$-graded algebra the 
growth of ordinary codimensions and graded codimensions
may differ. For example, if $A=M_k(F)\otimes F{\mathbb Z}_2$ with the canonical
${\mathbb Z}_2$-grading induced from group algebra $F{\mathbb Z}_2$, where $M_k(F)$ is
full $k\times k$ matrix algebra, then $exp(A)=k^2$ while $exp^{gr}(A)=2k^2$ (see
\cite{GZbook} for details). In the Lie case one can take $L=L_0\oplus L_1$ to be
a
two-dimensional metabelian algebra with $L_0=<e>, L_1=<f>$ and with  only one
non-trivial product $[e,f]=f$. Then $c_n(L)=n-1$ for all $n\ge 2$ hence $exp(L)=1$. On
the other hand $exp^{gr}(L)=2$.

\section{The main constructions and definitions}

Let $A$ be an arbitrary non-associative algebra over a field $F$ and let $F\{X\}$ be
an
absolutely free $F$-algebra with a
countable generating set $X$. A polynomial
$f=f(x_1,\ldots,x_n)$ is said to be an identity of $A$ if $f(a_1,\ldots,a_n)=0$ for
any $a_1,\ldots,a_n\in A$. The set of all identities of $L$ forms a T-ideal $Id(A)$
in $F\{X\}$, that is an ideal which is
stable under all endomorphisms of $F\{X\}$. Denote by
$P_n=P_n(x_1,\ldots,x_n)$ the subspace of all multilinear polynomials on
$x_1,\ldots,x_n$ in $F\{X\}$. Then $P_n\cap Id(A)$ is a subspace of all multilinear
identities of $A$ of degree $n$. In the case when ${\rm char}~F=0$,
the T-ideal $Id(A)$ 
is completely determined by the subspaces $\{P_n\cap Id(A)\}, n=1,2,\ldots~$.

For estimating how many identities an algebra $A$ can have 
one can define the so-called $n$-th
codimension of the identities of $A$ or, for shortness,
codimension  of $A$:
$$
c_n(A)=\dim \frac{P_n}{P_n\cap Id(A)},~
n=1,2,\ldots~.
$$

As it was mentioned above, the class of associative and non-associative algebras with
exponentially bounded sequence $\{c_n(A)\}$ is sufficiently wide. In the 
case when $c_n(A) < a^n$ for some real $a$, one can define the lower and the upper PI-exponents of $A$ as
follows:
$$
\underline{exp}(A)=\liminf_{n\to\infty} \sqrt[n]{c_n(A)},\quad
\overline{exp}(A)=\limsup_{n\to\infty} \sqrt[n]{c_n(A)}
$$

and the ordinary PI-exponent
\begin{equation}\label{e1}
exp(A)=\lim_{n\to \infty}\sqrt[n]{c_n(A)},
\end{equation}
provided that $\underline{exp}(A)=\overline{exp}(A)$.

For ${\mathbb Z}_2$-graded algebras one can also consider graded identities. Let $X$
and $Y$ be two infinite sets of variables and let $F\{X\cup Y\}$ be an
absolutely free
algebra generated by $X\cup Y$. If we suppose that all elements of $X$ are even and all
elements of $Y$ are odd, i.e. $\deg(x)=0, \deg(y)=1$ for any $x\in X,y\in Y$ then
$F\{X\cup Y\}$ can be naturally endowed by
a
${\mathbb Z}_2$-grading. A polynomial
$f=f(x_1,\ldots,x_m,y_1,\ldots, y_n)\in F\{X\cup Y\}$ is said to be a graded identity
of a superalgebra $A=A_0\oplus A_1$ if $f(a_1,\ldots,a_m,b_1,\ldots, b_n)=0$ for all
$a_1,\ldots,a_m\in A_0,b_1,\ldots, b_n\in A_1$. Fix $0\le k\le n$ and denote by
$P_{k,n-k}$ the subspace of $F\{X\cup Y\}$ spanned by all multilinear polynomials in
$x_1,\ldots,x_k \in X$, $y_1,\ldots, y_{n-k}\in Y$. Then $P_{k,n-k}\cap Id(A)$ is the
set of all multilinear polynomial identities of the superalgebra $A=A_0\oplus A_1$ in
$k$ even and $n-k$ odd variables.

One of the equivalent definitions of graded codimensions of $A$ is
$$
c_n^{gr}(A)=\sum_{k=0}^n {n\choose k} c_{k,n-k}(A),
$$
where
$$
c_{k,n-k}(A)=\dim\frac{P_{k,n-k}}{P_{k,n-k}\cap Id(A)}.
$$

Starting from a
${\mathbb Z}_2$-graded algebra of some class (Lie, Jordan alternative,
etc.) one can construct a
${\mathbb Z}_2$-graded algebra of different class using the
notion of the
Grassmann envelope. Grasmann envelopes play an exceptional role in
PI-theory. For example, any variety of associative algebras is generated by the
Grassmann envelope of some finite dimensional associative superalgebra \cite{Kem}. 
In Lie case any so-called special variety is generated by the Grassmann envelope of a
finitely generated Lie superalgera \cite{Va}.

We recall this construction for Lie and super Lie cases. Let $G$ be the Grassmann
algebra generated by $1$ and the infinite set $\{e_1,e_2,\ldots \}$ satisfying the
following relations: $e_i e_j=-e_j e_i, i,j=1,2,\ldots~$. It is known that $G$ has a
natural ${\mathbb Z}_2$-grading $G=G_0\oplus G_1$ where
$$
G_0=Span<e_{i_1}\cdots e_{i_n}|n=2k, k=0,1,\ldots>,
$$
$$
G_1=Span<e_{i_1}\cdots e_{i_n}|n=2k+1, k=0,1,\ldots>.
$$

Given a Lie algebra $L$ with ${\mathbb Z}_2$-grading $L=L_0\oplus L_1$, its Grassmann
envelope
$$
G(L)=L_0\otimes G_0\oplus L_1\otimes G_1 \subset L\otimes G
$$
is a Lie superalgebra. Vice versa, if $L=L_0\oplus L_1$ is a Lie superalgebra then
$G(L)$ is an ordinary Lie algebra with a
${\mathbb Z}_2$-grading.

\section{Cocharacters of Grassmann envelopes}

The main tool in studying codimensions asymptotics  is representation theory of
symmetric groups. We refer the reader to \cite{JK} for details. Symmetric group $S_n$
acts naturally on multilinear polynomials in $F\{X\}$ as
\begin{equation}\label{e2}
\sigma f(x_1,\ldots, x_n)= f(x_{\sigma(1)},\ldots, x_{\sigma(n)}).
\end{equation}
Hence $P_n$ is an $FS_n$-module and $P_n\cap Id(L)$ and also
$$
P_n(L)=\frac{P_n}{P_n\cap Id(L)}
$$
are $FS_n$-modules. $S_n$-character $\chi(P_n(L))$ is called $n$-th cocharacter of $L$
and we shall write
$$
\chi_n(L)=\chi(P_n(L)).
$$
Recall that any irreducible $FS_n$-module corresponds to a
partition $\lambda$ of $n$,
$\lambda\vdash n$, $\lambda=(\lambda_1,\ldots, \lambda_k)$, where
$\lambda_1\ge\ldots\ge \lambda_k$ are
positive integers and $\lambda_1+\cdots+ \lambda_k=n$. By the Maschke Theorem, any
finite dimensional $FS_n$-module $M$ decomposes into the direct sum of irreducible
components and hence its character $\chi(M)$ has a decomposition
$$
\chi(M)=\sum_{\lambda\vdash n} m_\lambda \chi_\lambda
$$
where $m_\lambda$ are non-negative integers. In particular, for the algebra $L$ we have
\begin{equation}\label{e3}
\chi(L)=\sum_{\lambda\vdash n} m_\lambda \chi_\lambda.
\end{equation}

Integers $m_\lambda$ in (\ref{e3}) are called multiplicities of $\chi_\lambda$ in
$\chi_n(L)$ and $d_\lambda=\deg \chi_\lambda=\chi_\lambda(1)$ are the dimensions of
corresponding irreducible representations. Therefore
\begin{equation}\label{e4}
c_n(L)=\dim P_n(L)=\sum_{\lambda\vdash n} m_\lambda d_\lambda.
\end{equation}

For any partition $\lambda=(\lambda_1,\ldots, \lambda_k)\vdash n$ one can construct
Young diagram $D_\lambda$ containing $\lambda_1$ boxes in the first row, $\lambda_2$
boxes in the second row and so on:
$$
D_\lambda=\; \begin{array}{|c|c|c|c|c|c|c|} \hline  &   &
\cdots &  &  & \cdots &   \\ \hline & & \cdots &  \\
\cline{1-4}\vdots  \\ \cline{1-1} \\ \cline{1-1} \end{array}
$$

Given integers $k,l,d \ge 0$, we define the partition
$$
h(k,l,d)=(\underbrace{l+d,\ldots, l+d}_{k},\underbrace{l,\ldots,l}_{d})
$$
of $n=kl+d(k+l)$. The Young diagram associated with $h(k,l,d)$ is hook shaped, and we define
$H(k,l)$, an infinite hook, as the union of all $D_\lambda$ with $\lambda=h(k,l,d)$, $d=1,2,\cdots~$. 
For shortness we will say that a partition $\lambda\vdash n$ lies in the hook $H(k,l)$,
$\lambda\in H(k,l)$, if $D_\lambda\subset H(k,l)$. In other words,
$\lambda\in H(k,l)$ if $\lambda=(\lambda_1,\cdots,\lambda_t)$ and $\lambda_{k+1}\le l$.
According to this definition we will say that the cocharacter of $L$ lies in the hook
$H(k,l)$ if $m_\lambda=0$ in (\ref{e3}) as soon as $\lambda\not\in H(k,l)$.

A particular case of $H(k,l)$ is an infinite strip $H(k,0)$. In this case $\lambda\in H(k,0)$
if $\lambda_{k+1}=0$.

The following fact is well-known and we state it without proof.

\begin{lemma}\label{l1}
Let $L$ be a finite dimensional algebra, $\dim L=d<\infty$. Then $\chi_n(L)$ lies in the hook
$H(d,0)$ for all $n\ge 1$.
\end{lemma}

\hfill $\Box$

Another important numerical invariant of the identities of $L$ is the colength $l_n(L)$.
By definition
\begin{equation}\label{e5}
l_n(L)=\sum_{\lambda\vdash n} m_\lambda
\end{equation}
where $m_\lambda$ are taken from (\ref{e3}). It easily follows from (\ref{e4}) and (\ref{e5}) that
\begin{equation}\label{e6}
\max\{d_\lambda|m_\lambda\ne 0\}\le c_n(L)\le l_n(L)\cdot\max\{d_\lambda|m_\lambda\ne 0\}.
\end{equation}

For studying graded identities of $L=L_0\oplus L_1$ we need to act separately on even
and odd variables. More precisely, the space $P_{k,n-k}= P_{k,n-k}(x_1,\ldots,
x_k,y_1,\ldots,y_{n-k})$ is an $S_k\times S_{n-k}$-module where symmetric groups $S_k$,
$S_{n-k}$ act on $x_1,\ldots, x_k$ and $y_1,\ldots,y_{n-k})$, respectively. Any
irreducible $S_k\times S_{n-k}$-module is a tensor product of $S_k$-module and an 
$S_{n-k}$-module and corresponds to the pair $\lambda, \mu$ of partitions,
$\lambda\vdash k, \mu\vdash n-k$. As before, the subspace $P_{n-k}\cap Id(L)$ is an 
$S_k\times S_{n-k}$-stable subspace and one can consider the quotient
$$
P_{k,n-k}(L)=\frac{P_{k,n-k}}{P_{k,n-k}\cap Id(L)}
$$
as an $S_k\times S_{n-k}$-module. Its $S_k\times S_{n-k}$-character $\chi_{k,n-k}(L)=
\chi(P_{k,n-k}(L))$ is decomposed into irreducible components.
\begin{equation}\label{e7}
\chi_{k,n-k}(L)=\sum_{{\lambda\vdash k \atop \mu\vdash n-k}} m_{\lambda,\mu} \chi_{\lambda,\mu}
\end{equation}
and we define the
$(k,n-k)$-colength of $L$ as
$$
l_{k,n-k}(L)=\sum_{{\lambda\vdash k \atop \mu\vdash n-k}} m_{\lambda,\mu}
$$
with $m_{\lambda,\mu}$ taken from (\ref{e7}).

First, we prove some relations between graded and non-graded numerical invariants. 
We begin by recalling
the correspondence between multilinear homogeneous polynomials in a
free
${\mathbb Z}_2$-graded Lie algebra and in a
free Lie superalgebra. Let $f=f(x_1,\ldots,x_k,
y_1,\ldots, y_m)$ be a non-associative polynomial multilinear on $x_1,\ldots,x_k$,
$y_1,\ldots, y_m$, where $x_1,\ldots,x_k$ are supposed to be even and $y_1,\ldots, y_m$
odd indeterminates. Then $f$ is a linear combination of monomials from $P_{k,m}$.
Let $M=M(x_1,\ldots,x_k,y_1,\ldots, y_m)$ be such a monomial. We fix positions of
$y_1,\ldots, y_m$ in $M$ and write $M$ for shortness in the following form
$$
M=X_0y_{\sigma(1)}X_1\cdots X_{m-1}y_{\sigma(m)}X_m
$$
where $X_0,\ldots, X_m$ are some words (possibly empty) consisting of
left and right brackets
and indeterminates $x_1,\ldots,x_k$. Now we define a monomial $\widetilde M$
on even indeterminates $x_1,\ldots,x_k$ and odd indeterminates $y_1,\ldots,y_m$ from
free Lie superalgebra as
$$
\widetilde M ={\rm sgn}(\sigma) X_0y_{\sigma(1)}X_1\cdots X_{m-1}y_{\sigma(m)}X_m.
$$
Extending this map $\,\widetilde \empty\,$ by linearity we obtain a linear isomorphism
$P_{k,m}\rightarrow P_{k,m}$ of two subspaces of
a
${\mathbb Z}_2$-graded free Lie
algebra and a
free Lie superalgebra, respectively. Although the monomials in $P_{k,m}$
are not linearly independent, it easily follows
from Jacobi and super-Jacobi identities 
that the map $\,\widetilde \empty\,$ 
is well-defined. Similarly, we can define the
inverse map from a free Lie superalgebra to a free ${\mathbb Z}_2$-graded Lie algebra.

Following the same argument as in the
associative case (see \cite[Lemma 3.4.7]{GZbook}) we
obtain for any ${\mathbb Z}_2$-graded Lie algebra $L$ and its Grassmann envelope
$G(L)=G_0\otimes L_0\oplus G_1\otimes L_1$ the following result.
\begin{lemma}\label{l2}
Let $f\in P_{k,m}$ be a multilinear polynomial in the free Lie algebra. Then
\begin{itemize}
\item
$f$ is a graded identity of $L$ if and only if $\widetilde f$ is a graded identity
of $G(L)$; and
\item
$\widetilde{\widetilde f} = f$.
\end{itemize}
\end{lemma}
\vskip -0.2in
\hfill $\Box$

The next lemma is an obvious generalization of Lemma \ref{l1}.
\begin{lemma}\label{l3}
Let $L=L_0\oplus L_1$ be a finite dimensional Lie algebra, $\dim L_0=k, \dim L_1=l$, and
let
$$
\chi_{q,n-q}(L)=\sum_{{\lambda\vdash q \atop \mu\vdash n-q}} m_{\lambda,\mu} \chi_{\lambda,\mu}
$$
be its $(q,n-q)$-graded cocharacter. If $m_{\lambda,\mu}\ne 0$ then $\lambda\in H(k,0)$
and $\mu\in H(l,0)$.
\end{lemma}
\hfill 
$\Box$

Using this remark we restrict the shape of the
graded cocharacter of the Grassmann envelope $G(L)$.

\begin{lemma}\label{l4}
Let $L=L_0\oplus L_1$ be a finite dimensional Lie algebra, $\dim L_0=k, \dim L_1=l$, and let $\widetilde L$ be its Grassmann envelope. If
\begin{equation}\label{e8}
\chi_{q,n-q}(\widetilde L)=\sum_{{\lambda\vdash q \atop \mu\vdash n-q}} m_{\lambda,\mu} \chi_{\lambda,\mu}
\end{equation}
and $m_{\lambda,\mu}\ne 0$ in (\ref{e8}) then $\lambda\in H(k,0)$ and $\mu\in H(0,l)$.
\end{lemma}
{\em Proof}. Suppose $m_{\lambda, \mu}\ne 0$ in (\ref{e8}) for some $\lambda\vdash q,
\mu \vdash n-q$. 
Then there exists a multilinear polynomial $g=g(x_1,\ldots,x_q,
y_1,\ldots,y_{n-q})$ such that
$$
f= e_{T_\lambda} e_{T_\mu} g(x_1,\ldots, y_{n-q})
$$
is not a graded identity of $\widetilde L$, where $e_{T_\lambda}\in FS_q,
e_{T_\mu}\in FS_{n-q}$ are essential idempotents generating minimal left ideals in
$FS_q, FS_{n-q}$, respectively. Inclusion $\lambda\in H(k,0)$ immediately follows by
Lemma \ref{l3} since $L$ and $G(L)$ have the same cocharacters on even indeterminates.
Since $e_{T_\lambda}$ and $e_{T_\mu}$ commute, applying Lemma 4.8.6 from \cite{GZbook}
we get
$$
\widetilde f= a e_{T_\lambda} g,
$$
where $a\in I_{\mu'}$. Here $\mu'$ is the conjugated to $\mu$ partition of $n-q$ and
$I_{\mu'}$ is the minimal two-sided ideal of $FS_{n-q}$ generated $e_{T_{\mu'}}$. That is,
$I_{\mu'}$ has the character $r\cdot\chi_{\mu'}$, where $r=d_{\mu'}=\deg \mu'$.

By Lemma \ref{l2}, $\widetilde f$ is not a graded identity of $G(\widetilde L)$. Since
$\widetilde{\widetilde h}=h$ for any $h\in P_{q,n-q}$, we see that $\widetilde f$ is not a graded
identity of $L$ and $\mu'\in H(l,0)$ by Lemma \ref{l3}. In other words, the number
of rows of Young diagram $D_{\mu'}$ does not exceed $l$. This number equals  the
number of columns of $D_\mu$ hence $\mu\in H(0,l)$ and we are done.
\hfill $\Box$

Using the previous lemma we restrict the shape of non-graded cocharacter of $G(L)$.

\begin{lemma}\label{l5}
Let $L=L_0\oplus L_1$ be a finite dimensional Lie algebra, $\dim L_0=k, \dim L_1=l$, and let
$$
\chi(\widetilde L)=\sum_{\lambda\vdash n} m_{\lambda} \chi_{\lambda}
$$
be
the
$n$-th (non-graded) cocharacter of $\widetilde L= G(L)$. Then $m_\lambda\ne 0$ 
only if $\lambda\in H(k,l)$.
\end{lemma}

{\em Proof}. Suppose $f\in P_n$ is not an identity of $\widetilde L$. Since $f$ is
multilinear we may assume that $f(x_1,\ldots, x_q, y_1,\ldots, y_{n-q})\in P_{q,n-q}$
is not an identity of $\widetilde L$ for some $0\le q \le n$. Moreover, we can consider only the case when a graded polynomial $f$ generates in $P_{q,n-q}$ an irreducible $S_q\times S_{n-q}$-submodule $M$ with the character $(\chi_\lambda,\chi_\mu)$, $\lambda \vdash q, \mu\vdash n-q$.

Now we lift 
the
$S_q\times S_{n-q}$-action up to an
$S_n$-action and consider a
decomposition
of $FS_n M$ into irreducible components:
$$
\chi(FS_n M)= \sum_{\nu\vdash n} m_\nu \chi_\nu.
$$
Since $\lambda$ lies in the hook $H(k,0)$, i.e. the
horizontal strip of 
height $k$ by Lemma \ref{l4} and $\mu$ lies in $H(0,l)$, the vertical strip of 
width $l$, it follows
from the
Littlewood-Richardson rule for induced representations (\cite[2.8.13]{JK}, see also
\cite[Theorem 2.3.9]{GZbook}) 
that $m_\nu=0$ as soon as $\nu\not\in H(k,l)$
and we have completed the proof.
\hfill $\Box$

\begin{lemma}\label{l6}
Let $G(L)=\widetilde L=\widetilde L_0\oplus \widetilde L_1$ be the Grassmann envelope
of a finite dimensional  Lie algebra $L=L_0\oplus L_1$ with $\dim L_0=k, \dim L_1=l$. Then its colength sequence $\{l_n( \widetilde L)\}$ is polynomially
bounded.
\end{lemma}

{\em Proof.}
We use the notation $\{z_1,z_2,\ldots\}$ for non-graded indeterminates here since
$\{x_1,x_2,\ldots\}$ were even variables in the previous statements.

Let
\begin{equation}\label{eqq1}
\chi(\widetilde L)=\sum_{\lambda\vdash n} m_\lambda\chi_\lambda
\end{equation}
be the 
$n$-th cocharacter of $\widetilde L$. By Lemma~\ref{l5} we have 
$\lambda\in H(k,l)$ as soon as $m_\lambda\ne 0$ in (\ref{eqq1}). Fix
$\lambda\vdash n$ with $m_\lambda =m\ne 0$ and consider
the
$FS_n$-submodule
\begin{equation}\label{eqq2}
W_1\oplus\cdots\oplus W_m\subseteq P_n(\widetilde L)
\end{equation}
with $\chi(W_i)=\chi_\lambda,$ for all $i=1,\ldots, m$. 

We shall prove that
\begin{equation}\label{eqq3}
m\le (k+l) 2^{2kl}n^{k^2+l^2}
\end{equation}
in (\ref{eqq2}). Denote by $\lambda'_1,\ldots, \lambda'_l$ the heights of
the
first
$l$ columns of 
the
Young diagram $D_\lambda$. 
Clearly, it suffices
to prove the
inequality (\ref{eqq3}) only for $\lambda$ with $\lambda_k>l$ and 
$\lambda'_l>k$. Otherwise $\lambda\in H(k',l')$ with $k'\le k, l'\le l$ and
$k'+l'<k+l$.

Denote
$$
\mu_1=\lambda'_1-k,\ldots,\mu_l=\lambda'_l-k.
$$
Then $\lambda_1+\cdots+\lambda_k+\mu_1+\cdots+\mu_l=n$.

It is well-known (see, for example, \cite {MZ5}) that one can  choose multilinear
$f_1\in W_1,\ldots, f_m\in W_m$ such that $FS_nf_1=W_1,\ldots, FS_nf_m=W_m$ 
and each $f_i, i=1,\ldots, m$,  is symmetric on $k$ sets of indeterminates of orders
$\lambda_1,\ldots, \lambda_k$ and is alternating on $l$ sets of orders
$\mu_1,\ldots, \mu_l$.

According to this decomposition into
symmetric and alternating sets we rename
$z_1,\ldots, z_n$ as follows
\begin{equation}\label{eqq4}
\{z_1,\ldots, z_n\}=\{z^1_1,\ldots,z^1_{\lambda_1},\ldots, z^k_1,\ldots,z^k_{\lambda_k},
\bar z^1_1,\ldots,\bar z^1_{\mu_1},\ldots, \bar z^l_1,\ldots,\bar z^l_{\mu_l}
\},
\end{equation}
where each $f_i$ is symmetric on any set $\{z^j_1,\ldots,z^j_{\lambda_j} \}$, 
$j=1,\ldots,k$, and is alternating on any set 
$\{\bar z^s_1,\ldots,\bar z^s_{\mu_s}\}$, $s=1,\ldots, l$.

We shall find  $\delta_1,\ldots, \delta_m\in F$ such that
$$
f=\delta_1f_1+\cdots + \delta_mf_m
$$
is an identity of $\widetilde L$ if (\ref{eqq3}) does not hold.
Note that for any $\delta_1,\ldots, \delta_m\in F$ a polynomial
$f$ is also symmetric on each subset $\{z^i_1,\ldots,z^i_{\lambda_i}\}, 
1\le i\le k$, and alternating on each subset $\{\bar z^s_1,\ldots,\bar 
z^s_{\mu_s}\}$, $s=1,\ldots, l$.

Let $E=\{e_1,\ldots, e_{k+l}\}$ be a homogeneous basis of $L$  with
$E_0=\{e_1,\ldots, e_k\}\subset  L_0$, $E_1=\{e_{k+1},\ldots, e_{k+l}\}\subset L_1$.
Then $f$ is an identity of $\widetilde L$  if and only if $\varphi(f)=0$ for any
evaluation $\varphi: Z\to \widetilde L$  such that $\varphi(z_i)=g_i\otimes a_i,
1\le i \le n,$ where $a_i$ is a basis element from $E$ and $g_i\in G$  has the
same parity as $a_i$ and $g_1\cdots g_n\ne 0$ in $G$. 

Note also that $\varphi(f)=0$
implies $\varphi'(f)=0$ for any evaluation $\varphi'$ such that $\varphi'(z_i)= 
g_i'\otimes a_i, 1\le i\le n$ provided that $g_1\cdots g_n\ne 0$.

Using these two remarks we shall find an upper bound for the number of evaluations for asking the question whether $f$
is
an identity of $\widetilde L$ or not.

Consider first one symmetric  subset $Z_1=\{z^1_1,\ldots,z^1_{\lambda_1} \}$. If
$\varphi(z^1_i)=g\otimes e, \varphi(z^1_j)=h\otimes e,$ for some $i\ne j$ with 
$e\in E_1,$ then $\varphi(f)=0,$ as follows from the symmetry on $Z_1$. Hence we
need to check only evaluations with at most $r\le l$ odd values 
$\varphi(z^1_{i_1})=g_1\otimes e_{t_1}, \ldots, \varphi(z^1_{i_r})=g_r\otimes e_{t_r}$, where $e_{t_1},\ldots, e_{t_r}\in E_1$ are distinct. Since $Z_1$ is the 
symmetric set of
variables, the result of evaluation $\varphi$ does not depend (up to the sign) 
on the choice of $i_1,\ldots, i_r$. Hence we have ${l\choose r}$ possibilities.

Given $0\le r\le l$, we estimate the number of evaluations of remaining
$\lambda_1-r$ variables in the even component of $\widetilde L$.
First, let $r=0$ and $\varphi(z^1_i)=g_i\otimes a_i, a_i\in E_0, 1\le i\le 
\lambda_1$. If $e_1$ appears in the row $(a_1,\ldots, a_{\lambda_1})$ exactly
$\alpha_1$ times, $e_2$ appears $\alpha_2$ times and so on,
then the result
of such substitution depends only on $\alpha_1,\ldots, \alpha_k$ since $f$ is 
symmetric on $Z_1$. Hence we have no more than $(\lambda_1+1)^k$ variants
since $0\le \alpha_1,\ldots, \alpha_k\le \lambda_1$. In particular, we need at most
$(n+1)^k$ evaluations if $r=0$.

Let now $r=1$. We can replace by odd element an arbitrary variable from $Z_1$ and
get (up to the sign) the same value $\varphi(f)$ since $f$ is symmetric on $Z_1$.
Suppose say, that
$\varphi(z^1_{\lambda_1})=h\otimes e, e\in E_1$, and
$\varphi(z^1_1)=g_1\otimes a_1,\ldots, \varphi(z^1_{\lambda_1-1}) =
g_{\lambda_1-1}\otimes a_{\lambda_1-1},$ 
where all $a_j$ are even. 
If
$\alpha_1,\ldots,\alpha_k$ are the same integers as in the
case $r=0$ then the result
of the substitution also depends only on $\alpha_1,\ldots,\alpha_k$. Hence for $r=1$
we have at most
$$
{l\choose 1}\lambda^k_1 \le {l\choose 1}(n+1)^k
$$
variants for $\varphi$ since $0\le \alpha_1,\ldots,\alpha_k\le \lambda_1-1$.

Similarly, for general $0\le r\le l$ we have at most
$$
{l\choose r}(\lambda_1+1-r)^k \le {l\choose r}(n+1)^k
$$
variants. Therefore for evaluating all variables from $Z_1$ it 
suffices 
$$
\sum_{r=0}^l {l\choose r}(n+1)^k = 2^l(n+1)^k
$$
substitutions and for all symmetric variables we need at most
$$
(2^l(n+1)^k)^k
$$
substitutions.

Now consider the alternating set $Z_1'=\{\bar z^1_1,\ldots, \bar z^1_{\mu_1}\}$.
If $\varphi(\bar z^1_i)=g\otimes e, \varphi(\bar z^1_j)=h\otimes e,$
for some
$i\ne j$ with the same $e\in E_0,$
then $\varphi(f)=0$, hence we can choose only
$0\le r\le k$ distinct basis elements $b_1,\ldots,b_r\in E_0$ for values of
$\bar z^1_{i_1}, \ldots, \bar z^1_{i_r}$ of the type $g_i\otimes b_i$. Up to the
sign,
the result of the
substitution does not depend on $i_1,\ldots, i_r$ and we have only 
${k\choose r}$ options.

Suppose now that all $\varphi(\bar z^1_i), 1\le i\le r$, are fixed even values. Let
$$
\varphi(\bar z^1_{r+1})=g_1\otimes b_1,\ldots,
\varphi(\bar z^1_{\mu_1})=g_{\mu_1-r}\otimes b_{\mu_1-r},\quad b_1\ldots,
 b_{\mu_1-r}\in E_1.
$$
Then (up to the sign) the result of $\varphi$ depends only on the number of entries
of $e_{k+1},\ldots, e_{k+l}$ into the row $(b_1,\ldots, b_{\mu_1-r})$. Hence we have
at most $(\mu_1-r+1)^l$ variants for substitution of odd variables. As in the 
symmetric case we have the following upper bound
$$
\sum_{r=0}^k {k\choose r}(n+1)^l = 2^k(n+1)^l 
$$
for one subset and $(2^k(n+1)^l)^l$ for all skew variables.

We have proved that one can find $T\le 2^{kl}(n+1)^{l^2+k^2}$ evaluations
$\varphi_1,\ldots,\varphi_T$ such that the relations
\begin{equation}\label{eqq5}
\varphi_1(f)=\cdots=\varphi_T(f)=0
\end{equation}
imply $\varphi(f)=0$ for any evaluation $\varphi$, that is $f$ is an identity
of $\widetilde L$. Recall that $f=\delta_1f_1+\cdots+\delta_mf_m$. Therefore for any 
evaluation $\varphi$ the equality $\varphi(f)=0$ can be viewed as a system of $k+l$ 
homogeneous linear equations in the algebra $\widetilde L$ on unknown coefficients
$\delta_1,\ldots,\delta_m$. If (\ref{eqq3}) does not hold then the system 
(\ref{eqq5}) has a non-trivial solution $\bar\delta_1,\ldots,\bar\delta_m$ and
$f=\bar\delta_1f_1+\cdots+\bar\delta_mf_m$ is an identity of $\widetilde L$, a
contradiction. 

We have proved the inequality (\ref{eqq3}). From this inequality it follows that
all multiplicities in (\ref{eqq1}) are bounded by $(k+l)2^{2kl}n^{k^2+l^2}$. 
Finally
note that the number of partitions $\lambda\in H(k,l)$ is bounded by $n^{k+l}$.
Hence
$$
l_n(\widetilde L)<(k+l)2^{2kl}n^{k^2+l^2+kl}
$$
and we have thus completed the proof.

\hfill $\Box$
\medskip

As a corollary of previous results we obtain the following:

\begin{proposition}\label{p1}
Let $L=L_0\oplus L_1$ be a finite dimensional ${\mathbb Z}_2$-graded Lie algebra
with $\dim L_0=k, \dim L_1=l$ and let $\widetilde L=G(L)$ be its Grassmann envelope.
Then there exist constants $\alpha,\beta\in {\mathbb R}$ such that
$$
c_n(\widetilde L)\le \alpha n^\beta(k+l)^n.
$$
In particular,
$$
\overline{exp}(\widetilde L)=\limsup_{n\to\infty} \sqrt[n]{c_n(\widetilde L)} \le k+l.
$$
\end{proposition}
{\em Proof}. By \cite[Lemma 6.2.5]{GZbook},
there exist constants $C$ and $r$ such that
$$
\sum_{\lambda\in H(k,l)} d_\lambda \le C n^r(k+l)^n
$$
for all $n=1,2,\ldots~$. In particular,
$$
\max\{d_\lambda|\lambda\vdash n, \lambda\in H(k,l)\} \le C n^r(k+l)^n.
$$
Now Lemma \ref{l6} and the inequality (\ref{e6}) complete the proof.
\hfill $\Box$

\section{Existence of PI-exponents}

\begin{proposition}\label{p2}
Let $L$ be a finite dimensional simple Lie algebra over an algebraically
closed field of characteristic zero with some ${\mathbb Z}_2$-grading,
$L=L_0\oplus L_1$, $\dim L_0=k, \dim L_1=l$. Let also $\widetilde L=G(L)$ be its
Grassmann envelope. Then there exist constants $\gamma>0,\delta\in {\mathbb R}$
such that
$$
c_n(\widetilde L)\ge \gamma n^\delta(k+l)^n.
$$
In particular,
$$
\underline{exp}(\widetilde L)=\liminf_{n\to\infty} \sqrt[n]{c_n(\widetilde L)} \ge k+l.
$$
\end{proposition}
{\em Proof}. Denote $d=k+l=\dim L$. By \cite[Theorem 12.1]{Razm},
for the adjoint
representation of $L$ there exists a multilinear asssociative polynomial
$h=h(u^1_1,\ldots, u^1_d,\ldots,$ $ u^m_1,\ldots, u^m_d)$ alternating on each subset of
indeterminates $\{u^i_1,\ldots, u^i_d\}$, $1\le i\le m$, such that under any evaluation
$\varphi: u^i_j\to ad~ b^i_j, b^i_j\in L$, the value $\varphi(h)$ is a scalar linear
transformation of $L$ and $\varphi(h)\ne 0$ for some $h$. It follows that for
any integer $t\ge 1$ there exists a multilinear Lie polynomial
$$
f_t=f_t(u^1_1,\ldots, u^1_d,\ldots, u^{mt}_1,\ldots, u^{mt}_d, w)
$$
alternating on each set $\{u^i_1,\ldots, u^i_d\}$, $1\le i\le mt$ such that
$\varphi(f_t)\ne 0$ for some evaluation $\varphi: \{u^1_1,\ldots,u^{mt}_d,w\}\to
L_0\cup L_1$. Since $f_t$ is multilinear and alternating on each set
$\{u^i_1,\ldots, u^i_d\}$ and  $d=\dim L_0+\dim L_1$ it follows that for any $t\ge 1$
we get a graded multilinear polynomial
$$
f_t=f_t(x^1_1,\ldots, x^1_k,\ldots, x^{mt}_1,\ldots, x^{mt}_k,
y^1_1,\ldots, y^1_l,\ldots, y^{mt}_1,\ldots, y^{mt}_l,w)
$$
which is not a graded identity of $L$ and it is alternating on each subset
$\{x^i_1,\ldots, x^i_k \}$ and on each subset $\{y^i_1,\ldots, y^i_l \}$,
$1\le i\le mt$, where $x^i_j$'s are even and $y^i_j$'s are odd variables. The latter
indeterminate $w$ can be taken of arbitrary parity, say, $w=x_0$ is even.

Consider an
$S_p\times S_q$-action on
$$
P_{p+1,q}=P_{p+1,q}(x_0,x^1_1,\ldots, x^{mt}_k, y^1_1,\ldots, y^{mt}_l),
$$
where $p=mtk, q=mtl$ and $S_p, S_q$ act on $\{x^i_j \}$, $\{y^i_j \}$, respectively.
It follows from Lemma \ref{l3} 
that the
$S_p\times S_q$-character of the
submodule generated
by $f$ in $P_{p+1,q}$ lies in the pair of of strips $H(k,0)$, $H(l,0)$, that is
$$
\chi(F[S_p\times F_q]f) = \sum_{{\lambda\vdash p\atop \mu\vdash q}} m_{\lambda,\mu} \chi_{\lambda,\mu}
$$
with $m_{\lambda,\mu}=0$, unless $\lambda \in H(k,0), \mu \in H(l,0)$. Hence $\lambda$
is a partition of $mtk$ with at most $k$ rows. On the other hand, $f$ depends on $mt$
alternating subsets of even indeterminates of order $k$ each. It is well-known that in this
case $m_{\lambda,\mu}=0$ if $\lambda=(\lambda_1, \lambda_2,\ldots )$ and $\lambda_1
\ge mt+1$. It follows that only rectangular partition
\begin{equation}\label{e10}
\lambda=(\underbrace{mt,\ldots, mt}_k)
\end{equation}
can appear in $F[S_p\times F_q]f$ with non-zero multiplicity. Similarly,
\begin{equation}\label{e11}
\mu=(\underbrace{mt,\ldots, mt}_l)
\end{equation}
if $m_{\lambda,\mu}\ne 0$. Hence we can assume that $f$ has the form
$$
f=e_{T_\lambda} e_{T_\mu} g(x^1_1,\ldots, y^{mt}_l,w)
$$
with $\lambda$ and $\mu$ of the types (\ref{e10}), (\ref{e11}), respectively.

By Lemma \ref{l2}, the polynomial $\widetilde f$ is not an identity of the Lie
superalgebra $\widetilde L=G(L)$ and by Lemma 4.8.6 from \cite{GZbook}, the
graded
polynomial $\widetilde f$ generates in $P_{p+1,q}(\widetilde L)$ an irreducible
$S_{p}\times S_q$-submodule with the character $(\chi_\lambda,\chi_{\mu'})$, where
$$
\mu'=(\underbrace{l,\ldots, l}_{mt})
$$
is conjugated to a
$\mu$ partition of $mtl$.

First we apply Littelwood-Richardson rule and induce this  
$S_{p}\times S_q$-module up to $S_n$-module. Then we induce the
obtained
$S_n$-module up to $S_{n+1}$-module, where $n=p+q=mt(k+l)$. It follows from 
the
Littelwood-Richardson rule that the induced $S_{n+1}$-module
can contain only simple submodule corresponding to
partitions $\nu\vdash n+1$ such that the Young diagram $D_\nu$ contains a
subdiagram $D_{\nu_0}$, where 
$$
\nu_0=h(k,l,t_0)=(\underbrace{l+t_0,\ldots, l+t_0}_{k},\, 
\underbrace{l,\ldots, l}_{t_0})
$$
is a finite hook with $t_0\ge l-k, mt-kl$. Since we are interested
with 
asymptotic of codimensions, we may assume that $mt-kl>l-k$ and then $t_0=mt-kl$.
In particular, $\nu_0$ is a partition of $n_0=(k+l)t_0+kl$. Then $n+1-n_0 = 
(k+l-1)kl+1$ and by \cite[Lemma 6.2.4]{GZbook}
$$
d_{\nu_0} \le d_\nu \le n^c d_{\nu_0}
$$
where $c=(k+l-1)kl+1$ and
$$
d_{h(k,l,t_0)}\simeq an_0^b(k+l)^{n_0} \quad {\rm if}~ n_0\to\infty
$$
for some constants $a,b$, by Lemma 6.2.5 from \cite{GZbook}. 
Here the relation
$f(n)\simeq g(n)$ means that $\lim_{n\to \infty}\frac{f(n)}{g(n)}=1$.
Since $c_{n+1}(\widetilde L) \ge d_\nu$ we get the inequality
\begin{equation}\label{e12}
c_{n+1}(\widetilde L) \ge \alpha(n+1)^\beta(k+l)^{n+1}
\end{equation}
for all $n=m(k+l)t$, $t=1,2,\ldots$ for some constants $\alpha>0$ and $\beta$.

Since Lie algebra $L$ is simple, the Grassmann envelope $\widetilde L$ is a
centerless
Lie superalgebra. It is not difficult to see that $c_{r+1}(\widetilde L) \ge
c_r(\widetilde L)$ in this case for all $r\ge 1$. Hence by (\ref{e12}) we have
$$
c_{n+j}(\widetilde L) \ge \alpha (n+1)^\beta (k+l)^{n+1}
$$
for any $1\le j \le m(k+l)$. Since $n=m(k+l)t$ one can find constants
$\gamma >0$ and $\delta$ such that
$$
c_{r}(\widetilde L) \ge \gamma r^\delta (k+l)^{r}
$$
for all positive integer $r$ and we have completed
the proof.
\hfill $\Box$

Theorem \ref{t1}
now
easily follows from Propositions \ref{p1} and \ref{p2}.
\medskip

{\em Proof of Theorem \ref{t2}}. First we obtain an upper bound for
$c_n^{gr}(\widetilde L)$:
$$
c_n^{gr}(\widetilde L)=\sum_{q=0}^n {n\choose q} c_{q,n-q}(\widetilde L),
$$
where
\begin{equation}\label{e13}
c_{q,n-q}(\widetilde L)=\sum_{\lambda \vdash q\atop \mu\vdash n-q}
m_{\lambda,\mu} d_{\lambda,\mu}
\end{equation}
and $d_{\lambda,\mu}=\deg \chi_{\lambda,\mu}=\deg \chi_{\lambda}\cdot \deg
\chi_{\mu}=d_\lambda d_\mu$. Moreover, $\lambda\in H(k,0)$, $\mu\in H(0,l)$ by
Lemma \ref{l4}. Applying Lemma 6.2.5 from \cite{GZbook}, we obtain
$$
\sum_{\lambda \in H(k,0)\atop \lambda\vdash q} d_\lambda \le
 Cn^r k^q,\quad
\sum_{\mu \in H(0,l)\atop \mu\vdash n-q}d_\mu \le Cn^r l^{n-q}
$$
for some constants $C,r$ and hence
\begin{equation}\label{e14}
\sum_{\lambda \in H(k,0), \lambda\vdash q \atop \mu \in H(0,l), \mu\vdash n-q}
d_\lambda d_\mu \le C^2 n^{2r}k^q l^{n-q}
\end{equation}

On the other hand, graded colength
$$
l_{q,n-q}(\widetilde L)=\sum_{\lambda\vdash q \atop \mu\vdash n-q} m_{\lambda,\mu}
$$
is not greater than non-graded colength $l_n(\widetilde L)$. Since $l_n(\widetilde L)$
is polynomially bounded by Lemma \ref{l6}, one can find a polynomial $\varphi(n)$
such that
\begin{equation}\label{e15}
m_{\lambda,\mu}\le \varphi(n)
\end{equation}
for any $m_{\lambda,\mu}$ in (\ref{e13}). It now follows from (\ref{e13}), (\ref{e14}) and
(\ref{e15})
that for $\psi(n)= C^2 n^{2r}\varphi(n)$ we have
\begin{equation}\label{e16}
c_n^{gr}(\widetilde L)\le \psi(n)\sum_{q=1}^n {n\choose q} k^q l^{n-q}=
\psi(n)(k+l)^n
\end{equation}
and we have obtained an upper bound for $c_n^{gr}((\widetilde L))$.

On the other hand, in \cite[Lemma 3.1]{BGR} it is proved that for
any associative $G$-graded algebra $A$, where $G$ is a finite group,
an ordinary $n$-th codimension is less than
or equal to the
graded $n$-th
codimension, for any $n$. Proof of this lemma does not use
associativity. Hence
\begin{equation}\label{e17}
c_n^{gr}(\widetilde L)\ge c_n(\widetilde L).
\end{equation}
Theorem \ref{t2}
now
follows from (\ref{e16}), (\ref{e17}) and Proposition
\ref{p2} and we have completed the proof.
\hfill $\Box$

\end{document}